\documentclass[a4paper,12pt]{article}
\usepackage{amsfonts,amsmath,amsxtra}

\newtheorem{theorem}{Theorem}
\newtheorem{corollary}{Corollary}
\newtheorem{lemma}{Lemma}

\newenvironment{definition}
{\smallskip\noindent{\bf Definition\/}:}{\smallskip\par}

\newenvironment{example}
{\smallskip\noindent{\bf Example\/}.}{\smallskip\par}
\newenvironment{remark}
{\smallskip\noindent{\bf Remark\/}.}{\smallskip\par}

\newenvironment{proof}
{\noindent{\bf Proof\/}.}{{ $\Box$}\smallskip\par}

\title{Seifert cohomology of trees}
\author{E. Gorsky}
\begin{document}
\maketitle
\begin{abstract}
To every tree we associate a filtered cochain complex. Its cohomology and the corresponding spectral sequence have clear combinatorial description. If a tree is the Dynkin diagram of a simple plane curve singularity, the graded Euler characteristic of this complex coincides with the Alexander polynomial of the link. In this case we also point the relation to the Heegard-Floer homology theory, constructed by P. Ozsvath and Z. Szabo.
\end{abstract}

\section{Introduction}

In the series of articles (\cite{os},\cite{os1}) P. Ozsvath and Z. Szabo constructed and developed a theory of  Heegard-Floer knot homology. These homology have a lot of unexpected and interesting properties -- for example, they provide the exact method of the calculation of the genus of a knot. The Euler characteristic of the Heegard-Floer homology coincides with the Alexander polynomial of a knot. 

The main objects of our study will be algebraic knots and links, i. e. intersections of germs of curve singularities in $(\mathbb{C}^2,0)$ with a small sphere centered at the origin.
It turns out (\cite{os2}), that for algebraic knots the ranks of Heegard-Floer homology are completely determined by Alexander polynomial and can be deduced from it in a purely combinatorial procedure. 

Therefore one can ask if there is a homology theory based on the local topological information about a singularity (e. g. its Dynkin diagram)  coinciding with the Heegard-Floer homologies of its link. In this article we construct a combinatorial bigraded complex
whose homology coincide with the Heegard-Floer homology of a link for the simple irreducible plane curve singularities: $A_{2n}, E_6$ and $E_8$.

To every tree we associate a filtered cochain complex, which we call Seifert complex. Its cohomology (which we call Seifert cohomology) and the corresponding spectral sequence have clear combinatorial description. If a tree is the Dynkin diagram of a  plane curve singularity, the Euler characteristic of its Seifert cohomology coincides with the Alexander polynomial of the corresponding link (Theorem 2). We prove that, as for Heegard-Floer theory,  an analogue of the Poincare duality (Lemma 2) holds for Seifert cohomology, and the spectral sequence also converges to one-dimensional space (Theorem 1). We  compute the Seifert cohomology for several classes of trees, among then for all simply-laced Dynkin diagrams. 
In the last case we compare the answer with the Heegard-Floer homology of the corresponding knot.

\bigskip

The author is grateful to S. Gusein-Zade, V. Rubtsov, A. Gorsky, M. Kazaryan, S. Shadrin, G. Gusev, A. Kustarev and M. Bershtein for useful discussions and comments.
 
\section{Seifert complex}

Let T be an arbitrary tree. 

\begin{definition}
A {\it configuration} is a way of marking of some vertices of T by pluses and minuses and making some edges of T red such that every vertex of T is either marked by a plus or minus or it is the end of {\it exactly one} red edge.
For  a configuration $C$ we define the grading $Q(C)$ as the number of minuses and $E(C)$ as the sum of the number of minuses and the number of red edges. 
\end{definition}

\bigskip
 
Consider the vector space $SC(T)$ over $\mathbb{F}_2$ with the basis labeled by all possible configurations. The integer-valued functions $Q$ and $E$ make this space bigraded: we'll denote by $SC^{k}(T,n)$ the subspace spanned by configurations $C$ with $Q(C)=k$ and $E(C)=n$.

\medskip
 
Let us define the pair of differentials $D$ and $d$ on $SC(T)$ by the equations:

$$D(\odot\Longleftrightarrow\odot)=(\oplus \longleftrightarrow \ominus) +(\ominus \longleftrightarrow\oplus),\,\,\, d(\oplus)=\ominus.$$

These maps are extended to the space $SC(T)$ using the Leibnitz rule, and, thanks to our ground field $\mathbb{F}_2$,
we do not need to care about signs.

\begin{example}
Let $$C=(\oplus \longleftrightarrow \odot \Longleftrightarrow \odot \longleftrightarrow \ominus \longleftrightarrow \odot \Longleftrightarrow \odot \longleftrightarrow \oplus), $$
then 
$$D(C)=(\oplus \longleftrightarrow \oplus \longleftrightarrow \ominus \longleftrightarrow \ominus \longleftrightarrow \odot \Longleftrightarrow \odot \longleftrightarrow \oplus)+$$
$$(\oplus \longleftrightarrow \ominus \longleftrightarrow \oplus \longleftrightarrow \ominus \longleftrightarrow \odot \Longleftrightarrow \odot \longleftrightarrow \oplus)+$$
$$(\oplus \longleftrightarrow \odot \Longleftrightarrow \odot \longleftrightarrow \ominus \longleftrightarrow \oplus \longleftrightarrow \ominus \longleftrightarrow \oplus)+$$
$$(\oplus \longleftrightarrow \odot \Longleftrightarrow \odot \longleftrightarrow \ominus \longleftrightarrow \ominus \longleftrightarrow \oplus \longleftrightarrow \oplus);$$
$$d(C)=(\ominus \longleftrightarrow \odot \Longleftrightarrow \odot \longleftrightarrow \ominus \longleftrightarrow \odot \Longleftrightarrow \odot \longleftrightarrow \oplus)+$$
$$(\oplus \longleftrightarrow \odot \Longleftrightarrow \odot \longleftrightarrow \ominus \longleftrightarrow \odot \Longleftrightarrow \odot \longleftrightarrow \ominus).$$
\end{example}

\medskip

It is clear that $D$ and $d$ are differentials and commute, so $D+d$ is a differential too. 
It is also clear that $D$ changes the $(Q,E)$ bigrading by $(1,0)$ and $d$ changes it by $(1,1)$

\medskip

\begin{lemma}
The dimension of the cohomology of $SC(T)$ with respect to $d$ is less or equal to 1.
More precisely, if $T$ admits a decomposition in non-intersecting pairs of connected vertices, the above dimension is 1, otherwise it vanishes.     
\end{lemma} 

\begin{proof}
Let us fix the set of red edges and encounter all possible markings in free vertices. They form a binary cube, and $d$ is a natural differential on this cube. Therefore the complex $(SC(T),d)$ is a direct sum of complexes which are acyclic if and only if there are vertices that do not belong to any red edge. The cohomology of $(SC(T),d)$ is generated by decompositions of the whole $T$ into non-intersecting red edges.

It remains to remark that there cannot be more than one such decomposition: in any tree we can find a vertex of degree 1. 
An edge to it should be red, and by the induction we can reconstruct the configuration of red edges, if it's possible.
\end{proof}

\bigskip

The second differential has more interesting cohomology.

\begin{definition}
By the Seifert cohomology $SH^{\bullet}(T)$ of a tree $T$ we mean the cohomology of the complex $(SC^{\bullet}(T),D)$.
Since $D$ preserves the $E$-grading, we can decompose it into a sum
$$SH^{\bullet}(T)=\oplus_{n}SH^{\bullet}(T,n),$$
where $SH^{\bullet}(T,n)$ denotes the cohomology of the subcomplex  $(SC^{\bullet}(T,n),D)$.
\end{definition}

\bigskip

\begin{theorem}
There exists a natural spectral sequence starting from $SH^{\bullet}(T)$ and converging to a zero-, or one-dimensional $E_{\infty}$ term. Its differential $d_1$ is induced be $d$, and  higher differentials $d_n$ increase $E$-grading by $n$ and $Q$-grading by 1.
\end{theorem}

\begin{proof}
Since $d$ and $D$ commute, we can consider $SC^{\bullet}(T)$ as a bicomplex with gradings $E$ and $E-Q$. With this bicomplex we can associate the pair of spectral sequences both converging to the cohomology of the total complex
$$(SC^{\bullet}(T),d+D)$$ (with grading $Q$):
first has cohomology of $(SC^{\bullet}(T),d)$ as the $E_1$ term, and second has cohomology of $(SC^{\bullet}(T),D)$ as the $E_1$ term.
Lemma 1 says that the total dimension of $E_1$ term of the first spectral sequence is less or equal to 1, so it coincides with the $E_{\infty}$ term. The second spectral sequence starts from $SH^{\bullet}(T)$ and converges to the same $E_{\infty}.$  

It remains to remark that $d_n$ should increase $E$ by $n$ and $E-Q$ by $n-1$, so $Q$ is increased by 1.
\end{proof}

\bigskip

\begin{lemma}
The following natural duality holds:
$$SH^{k}(T,n)=SH^{V-2n+k}(T,V-n),$$
where $V$ is the number of vertices in $T$.
\end{lemma}

\begin{proof}
Consider the involution $*$ on $SC(T)$ which changes all pluses by minuses and vice versa. It commutes with $D$, so it gives an isomorphism of the corresponding Seifert cohomology. Now, if $E(C)=E$, $Q(C)=Q$, then $C$ has $Q$ minuses, $E-Q$ red edges, and $V-2E+Q$ pluses, so $$Q(*C)=V-2E+Q,\,\,\,\, E(*C)=V-E.$$ 
\end{proof}

\bigskip

Let us turn to the Euler characteristic of Seifert cohomology.

\begin{definition}
Let us enumerate the vertices of $T$ in an arbitrary way. Consider an upper-triangular matrix $S=(s_{ij})$ 
defined as follows:
$$s_{ij}=
\begin{cases}
-1, \mbox{\rm if}\,\,\,\,\, i=j\\
1, \mbox{\rm if}\,\,\,\, i<j\,\,\,\, \mbox{\rm and vertices $i$ and $j$ neighbor}\\
0, \mbox{\rm otherwise} 
\end{cases}
$$

We call $S$ the Seifert matrix of the tree T with the given enumeration of vertices.
\end{definition} 

\begin{theorem}
The generation function for the Euler characteristics of Seifert homologies equals to
$$\Delta(t)=\sum_{k,n}t^{n}\chi(SH^{\bullet}(T,n))=\det (tS-S^{T}).$$
\end{theorem}

\begin{proof}
First, let us note that $\Delta(t)$ is a sum over all configurations
$$\Delta(t)=\sum_{C}(-1)^{Q(C)}t^{E(C)}.$$
If we fix a set $R$ of red edges, and change arbitrary signs at $V-2|R|$ free vertices, we get the equation
$$\Delta(t)=\sum_{R}t^{|R|}\sum_{k=0}^{V-2|R|}{V-2|R|\choose k}(-1)^{k}t^{k}=\sum_{R}t^{|R|}(1-t)^{V-2|R|}.$$

On the other hand, consider the matrix $A=tS-S^{T}=(a_{ij})$. It has $(1-t)$ on the diagonal, $t$ above the diagonal, corresponding to edges, and $(-1)$ below the diagonal, also corresponding to edges. Consider an arbitrary permutation $\sigma$. If it contains a cycle $(i_1,\ldots,i_l)$ with $l>2$, then $a_{i_1i_2}\cdot\ldots\cdot a_{i_{l-1}i_{l}}\cdot a_{i_li_1}=0$, since otherwise the graph T should contain edges $i_1i_2,\ldots,i_{l-1}i_l,i_li_1$, and it's not a tree.

Therefore the determinant of $A$ equals to $\sum_{\sigma}(-1)^{\sigma}a_{1\sigma(1)}\ldots a_{V\sigma(V)}$, where the sum is taken over all products of non-intersecting transpositions along the edges of $T$. Such a permutation $\sigma$ is equivalent to some set of red edges $R$. The sign of $\sigma$ equals to $(-1)^{|R|}$, so
$$\det A=\sum_{R}(-1)^{|R|}(-t)^{|R|}(1-t)^{V-2|R|}=\Delta(t).$$ 
\end{proof}

\begin{corollary}
The determinant $\Delta_{T}(t)=\det(tS-S^{T})$ does not depend of the enumeration of vertices of T. We'll refer to it
as to the {\it Alexander polynomial} of the tree T. If we denote $M=(S^{T})^{-1}S$, then
$$\Delta_{T}(t)=det(E-tM).$$
\end{corollary}

\begin{lemma}
$$\Delta_{T}(1)=\dim E_{\infty}$$
\end{lemma}

\begin{proof}
Theorem 1 says that the dimension of the $E_{\infty}$ term of the spectral sequence is always 0 or 1.
We can reconstruct it from the Euler characteristic: this dimension is equal to $\Delta_{T}(1)$, since
the Euler characteristic of the $E_{\infty}$ term is nothing but $\sum_{C}(-1)^{Q(C)}=\Delta_{T}(1)$.
\end{proof}

\bigskip

The next two lemmas give the recursive description of the Alexander polynomial and Seifert cohomology with respect to forgetting of a degree 1 vertex of a tree.

\begin{lemma}
\label{hang}
Let a vertex $v$ of a tree $T$ has degree 1, and its neighbor is a vertex $w$. Let $T_1=T\setminus \{v\}$,
$T_2=T\setminus \{v,w\}.$ Then
$$\Delta_T(t)=(1-t)\Delta_{T_1}+t\Delta_{T_2}.$$
\end{lemma}

\begin{proof}
The edge $(vw)$ can be marked or not. If it is not marked, we choose arbitrarily the sign at $v$ and a configuration in $T_1$. If the edge is marked, we choose a configuration in $T_2$.
\end{proof}

The last lemma has a natural categorification.

\begin{lemma}
In the notation of Lemma \ref{hang}, the following short exact sequence of complexes exists for every $n$:
\begin{equation}
\label{recur}
0\rightarrow SC^{k}(T_1,n)\oplus SC^{k-1}(T_1,n-1)\rightarrow SC^{k}(T,n)\rightarrow SC^{k}(T_2,n-1)\rightarrow 0.
\end{equation}
\end{lemma}

\begin{proof}
The proof is analogous: we got a subcomplex $$SC^{k}(T_1,n)\oplus SC^{k-1}(T_1,n-1)$$ in the complex $SC^{k}(T,n)$, and the quotient is isomorphic to $SC^{k}(T_2,n-1)$.
\end{proof}

\section{Examples}

Let us turn to some examples of computation of Seifert cohomology. We will pack it into the Poincare polynomial
$$P_T(q,t)=\sum_{k,n}q^{k}t^{n}\dim SH^{k}(T,n).$$

\begin{example}
Let us calculate the part of Seifert cohomology of arbitrary tree T without red edges ("zero level" with $E=Q$).  From the definition of $D$ we see that the fragments $(\oplus \longleftrightarrow \ominus)$ and $(\ominus \longleftrightarrow\oplus)$ are cohomological,
so the unique invariant of the configuration is the total number of minuses (i.e. $E=Q$). 
Therefore for any tree $T$ we have
$$SH^{k}(T,k)=
\begin{cases}
 \mathbb{F}, 0\le k\le |T|\\
 0, \mbox{\rm otherwise}\\
\end{cases}
$$
\end{example}
 
\begin{example}
Let us compute the Seifert cohomology of the Dynkin diagram of type $A_n$. Let us prove that they are located only at "zero level" by induction on $n$ with the help of the exact sequence (\ref{recur}).
We have the long exact sequence of Seifert cohomology:

$$\ldots\rightarrow SH^{i-1}(A_{n-2},k-1)\rightarrow SH^{i}(A_{n-1},k)\oplus SH^{i-1}(A_{n-1},k-1)\rightarrow$$ $$\rightarrow SH^{i}(A_{n},k)\rightarrow SH^{i}(A_{n-2},k-1)\rightarrow\ldots$$

Since we know, that $SH^{j}(A_{l},k)$ is nonzero only for $0\le j=k\le l$ for all $l<n$, then for fixed $k$ the 
only nonzero part in the above exact sequence is

$$0\rightarrow SH^{k-1}(A_{n},k)\rightarrow SH^{k-1}(A_{n-2},k-1)\rightarrow$$
$$\rightarrow SH^{k}(A_{n-1},k)\oplus SH^{k-1}(A_{n-1},k-1)\rightarrow SH^{k}(A_n,k)\rightarrow 0.$$

For $k\le n-1$ we have $$\dim SH^{k-1}(A_{n-2},k-1)=\dim SH^{k}(A_{n-1},k)=\dim SH^{k-1}(A_{n-1},k-1)=1$$ by the
induction assumption,and $\dim SH^{k}(A_n,k)=1$ by the previous example. Hence we get $SH^{k-1}(A_{n},k)=0$.
For $k=n$ we have $$\dim SH^{k-1}(A_{n-2},k-1)=\dim SH^{k}(A_{n-1},k)=0,$$ $$\dim SH^{k-1}(A_{n-1},k-1)=\dim SH^{k}(A_n,k)=1,$$ so $SH^{k-1}(A_{n},k)=0$.

Therefore the Poincare polynomial of Seifert cohomology of $A_n$ diagram equals to
$$P_{A_n}(q,t)=1+qt+q^2t^2+\ldots+q^nt^n.$$
\end{example}

\begin{example}
Let us compute the Seifert cohomology of the Dynkin diagram  of type $D_4$. Among the configurations with $E=2$ there are 6  without red edges and 6  with one red edge. One can check that $D$ has a one-dimensional kernel on $SC^{\bullet}(D_4,2)$, generated by the sum of all configurations with one red edge, one plus and one minus. We'll refer to this sum as to {\it turbine}. 
The Poincare polynomial equals to
$$P_{D_4}(q,t)=(1+qt+q^2t^2+q^3t^3+q^4t^4)+qt^2.$$ 
\end{example}

\begin{example}
As a generalization of the previous example, let us compute the Seifert cohomology of the Dynkin diagram of type $D_n$. On the zero level the cohomology are known, the others have one red edge and can be described as follows (with the help of the exact sequence (\ref{recur})):
we take the turbine in the $D_4$ subdiagram, and in the remaining  $A_{n-4}$ subdiagram we take arbitrary generator on zero level. The Poincare polynomial equals to
$$P_{D_n}(q,t)=P_{A_n}(q,t)+qt^2P_{A_{n-4}}(q,t).$$
The Alexander polynomial equals to
$$\Delta(D_n)=\Delta(A_n)-t^2\Delta(A_{n-4})=1-t+(-1)^{n-1}t^{n-1}+(-1)^nt^n.$$
\end{example}

\begin{example}
Consider the Dynkin diagrams of types $E_6, E_7, E_8$. If we throw the $D_4$ subdiagram from $E_n$, we get the disconnected sum of  $A_1$ and $A_{n-5}$. Therefore
$$P_{E_n}(q,t)=P_{A_n}(q,t)+qt^2P_{A_1}(q,t)P_{A_{n-5}}(q,t)=P_{A_n}(q,t)+qt^2(1+qt)P_{A_{n-5}}(q,t).$$
\end{example}

\begin{example}
Let us compute the Seifert cohomology of the graph $S_n$ with one vertex of degree $n-1$ and others of degree 1.
Among all configurations with a given value of $E$ there are ${n\choose E}$ on zero level and $(n-1){n-2\choose E-1}$ 
on the with one red edge ("first level").
On the zero level cohomology are one-dimensional, so on the first level the cohomology dimension equals to
$(n-1){n-2\choose E-1}-{n\choose E}+1$. The Poincare polynomial equals to
$$P_{S_n}(q,t)=P_{A_n}(q,t)+\sum_{E=2}^{n-2}[(n-1){n-2\choose E-1}-{n\choose E}+1]q^{E-1}t^E=$$
$$P_{A_n}(q,t)+(n-1)t(1+qt)^{n-2}-q^{-1}[(1+qt)^{n}-P_{A_{n}}(q,t)]=$$
$$q^{-1}[(1+q)P_{A_n}-(1+qt)^{n-2}(1-(n-3)qt+q^2t^2)].$$
\end{example}

\bigskip

Now let us turn to the calculation of the spectral sequences.  
 
\begin{example}
Consider the $A_{n}$ diagram. By the Leibnitz's formula the generator with $k$ minuses,
is mapped by $d=d_1$ to the generator with $k+1$ minuses, with the coefficient $(n-k)$.
Since we work over $\mathbb{F}_2$, we have $d_1([n-2k])=0$, and $d_1([n-2k+1])=[2k+2]$.
Therefore $E_2$ term in the spectral sequence vanishes for $n$ odd and has one generator $[0]$ for even $n$,
so $E_2=E_{\infty}$. 
\end{example}

\begin{example}
For the $D_4$ diagram the spectral sequence is more interesting. The $E_2$ term is generated by $[0]$ and the turbine. 
We know that $E_{\infty}=0$, so the differential $d_2$ maps $[0]$ to the turbine.

This can be checked more explicitly. If $C_0$ is a configuration with all pluses, $d(C_0)$ is a sum of 4 configurations with one minus, $D^{-1}(d(C_0))$ is a sum of 3 configurations with one red edge and two minuses,
so $$d_2(C_0)=d(D^{-1}(d(C_0)))$$ equals to the turbine. 
\end{example}

\section{Seifert form}

A Seifert surface of a knot is an oriented surface $\Sigma$ in $S^3$ whose border is a given knot. To a Seifert surface one can naturally associate the bilinear form on  $H^{1}(\Sigma)$ defined in the following way: if $x$ and $y$ are two cycles on $\Sigma$, let $x^{+}$ be the shift of $x$ along the small positive normal vector field to $\Sigma$.
Define $<x,y>$ as the linking number of $x^{+}$ and $y$. The matrix $S$ of the bilinear form $<\cdot,\cdot>$ is called the Seifert matrix, and the Alexander polynomial  of the knot can be expressed via the Seifert matrix as
$$\Delta(t)=Det(tS-S^{T}).$$

By an algebraic knot we mean the intersection of a germ of plane curve $\{f=0\}$ in $(\mathbb{C}^2,0)$ with a small
sphere centered at the origin. 
For the algebraic knots one can choose the Milnor fiber $\{f=\varepsilon\}$ as a natural Seifert surface.
In its homologies one can define the class of distinguished basises of vanishing cycles (\cite{book}).

\begin{theorem}(\cite{book})
In a distinguished basis the following facts hold:

1. The Seifert matrix $S$ is upper-triangular with the numbers $(-1)$ on the diagonal.

2. The matrix $S-S^{T}$ coincides with the intersection matrix in $H^{1}(V_{\varepsilon})$ in this basis.

3. The matrix $M=(S^{T})^{-1}S$ coincides with the  matrix of the classical monodromy in this basis.
\end{theorem}

The expression for the Alexander polynomial follows from the proposition 3 since
$$\Delta(t)=\det (E-tM).$$

Moreover, it is well known that for the simple singularities one can provide the distinguished bases with the intersection matrix corresponding to the Dynkin diagram of the same name as the singularity ($A_n,D_n,E_6,E_7,E_8$). Therefore the Seifert matrix in the above meaning coincide with the matrix of the Seifert form in the distinguished basis.

\begin{lemma}
If a tree T is a Dynkin diagram of an irreducible singularity, then the spectral sequence from Theorem 1 converges to one-dimensional space.

If a singularity is reducible (i. e. we get a link), its limit term vanishes.
\end{lemma}

\begin{proof}
It follows from Theorem 2, that the dimension of $E_{\infty}$ term of this spectral sequence equals to $\Delta(1)$.
Now $\Delta(1)=\det(S-S^{T})$, what equals to the determinant of the intersection matrix. It is equal to 1, if the border of a Seifert surface has one component, and to 0 otherwise.
\end{proof}

\section{Heegard-Floer homologies} 

In the series of articles P. Ozsvath and Z. Szabo (\cite{os},\cite{os1}) constructed a categorification of the Alexander polynomial with the following properties. For every knot $K$ there exists a filtered complex $CF^{-}(K)$ of free $\mathbb{Z}[U]$ -- modules such that:

1. The operator $U$ increases the homological grading by 2, and the filtration level by 1.

2. $$H^{\bullet}(CF^{-}(K))\simeq \mathbb{Z}[U], H^{\bullet}(CF^{-}(K)/UCF^{-}(K))=\mathbb{Z}.$$

Let $$\widehat{CF}(K)=CF^{-}(K)/UCF^{-}(K), HF^{-}(K,n)=H^{\bullet}(CF^{-}(K,n)/CF^{-}(K,n+1)),$$ $$\widehat{HF}(K,n)=H^{\bullet}(\widehat{CF}(K,n)/\widehat{CF}(K,n+1)).$$

3. $$\sum_{n=0}^{\infty}\chi(CF^{-}(K,n)/CF^{-}(K,n+1))t^n={\Delta_{K}(t)\over 1-t},$$ $$\sum_{n=0}^{\infty}\chi(\widehat{CF}(K,n)/\widehat{CF}(K,n+1))t^n=\Delta_{K}(t),$$
where $\Delta_K(t)$ is the Alexander polynomial of the knot $K$.

\begin{remark}
From the equation 2 there exist two natural spectral sequences. One starts from $\oplus_n \widehat{HF}^{\bullet}(K,n)$ and converges to $\mathbb{Z}$, second starts from $\oplus HF^{-}_n(K,n)$ and converges to 
$\mathbb{Z}[U]$.
\end{remark}

\medskip
 
For the algebraic knots the cohomology of the complexes 
$CF^{-}(K)$ and $\widehat{CF}(K)$ has the following explicit description (\cite{os2}). Following \cite{my}, 
we slightly modify for simplicity the gradings in the Heegard-Floer homologies. 

Let $${\Delta(t)\over 1-t}=1+t^{\alpha_1}+t^{\alpha_2}+\ldots,$$
consider $$P_{g}(t,q)=1+qt^{\alpha_1}+q^2t^{\alpha_2}+\ldots, \Delta_{g}(t,q)=(1-qt)P_g(t,q).$$
Let us change in $P_g$ the expression $q^k$ by $u^{2k}$, and in $\Delta_{g}$ we change $q^{k}$ by $u^{2k}$ and $-q^{k}$ by $u^{2k-1}$. We'll get the graded Poincare polynomial for $HF^{-}(K)$ and
$\widehat{HF}(K)$ respectively.

\begin{example}
Consider the singularity of type $E_6$, corresponding to the function $f=x^3-y^4$. 
Its link is the $(3,4)$ torus knot.
One can check that
$${\Delta(t)\over 1-t}=1+t^3+t^4+{t^6\over 1-t},$$
so
$$P_{g}(t,q)=1+qt^3+q^2t^4+{q^3t^6\over 1-qt},$$
$$\Delta_{g}(t,q)=1-qt+qt^3-q^3t^5+q^3t^6,$$
and the Poincare polynomial for the Heegard-Floer homology equals to
$$\widehat{HF}(t,u)=1+ut+u^2t^3+u^5t^5+u^6t^6.$$
\end{example}

\bigskip

There is a natural question of the description of Heegard-Floer homology of an algebraic knot in terms of the topological invariants of singularity. In \cite{my} a quite formal method was suggested, based on the results of A. Campillo, F. Delgado and S. Gusein-Zade   (e. g. \cite{cdg},\cite{cdg3}) on the motivic integration over the space of functions. It turns out, that for simple irreducible singularities one can naturally extract the Heegard-Floer homology from the above construction of Seifert cohomology .   

Consider the "degeneration" operators $K_n([m])=d_n([m-n])$. The operator $K_n$ does not change $E$-grading, and increases $Q$-grading by 1.
 
\begin{example}
Let us describe the action of $K_2$ in the Seifert cohomology of the Dynkin diagram of type $E_6$. These cohomology have 11 generators,
7 of them have a form $[k], 0\le k\le 6$, their $Q$ and $E$-gradings both equal to $k$. Let $[\tau]$ be the turbine 
corresponding to the $D_4$ subdiagram. Seifert cohomology has  4 other generators:
$[0]\otimes[\tau]\otimes[0],[1]\otimes[\tau]\otimes[0],[0]\otimes[\tau]\otimes[1],[1]\otimes[\tau]\otimes[1]$
with $Q$-gradings 1, 2, 2, 3 and $E$-gradings  2, 3, 3, 4. 

If we take the cohomology with respect to  $K_2$, we get the generators $[0], [1], [5], [6]$ and $[1]\otimes[\tau]\otimes[0]+[0]\otimes[\tau]\otimes[1]$ with $Q$-gradings 0,1,5,6,2 and $E$-gradings 0,1,5,6,3. This is exactly the same as the Heegard-Floer homologies obtained in the previous example. 
\end{example} 

\begin{example}
Let us make the analogous computation for the singularity $E_8$. After taking the cohomology of $K_2$ we get the generators  $$[0],[1],[1]\otimes[\tau]\otimes[0]+[0]\otimes[\tau]\otimes[1],[1]\otimes[\tau]\otimes[1]+[0]\otimes[\tau]\otimes[2],
[1]\otimes[\tau]\otimes[2]+[0]\otimes[\tau]\otimes[3],[7],[8]$$ of $Q$-gradings 0,1,2,3,4,7,8 and $E$-gradings 0,1,3,4,5,7,8.
This also coincides with the Heegard-Floer homology of the $(3,5)$ torus knot.
\end{example}

\section{Remarks and discussions}

1. In \cite{vc} S. Cecotti and C. Vafa gave the physical interpretation of the Seifert form and the Dynkin diagram.
The vertices of a Dynkin diagram are interpreted as "vacuum states" and edges correspond to "solitons". Extending this model, we can naturally identify red edges with the soliton-antisoliton pairs (since  they are related to transpositions in the proof of Theorem 2), and minuses may correspond to the ionizations of the vacuum states.
Now the grading $E$ may be interpreted as the energy of such a composite state of a physical system, and $Q$ can be naturally identified with its charge (since soliton-antisoliton pair is free of charge but has some energy).
Now $D$ is the operator of the decay of a soliton-antisoliton pair, and $d$ is an ionization operator.
The operator $D+d$, appearing in Theorem 1 as a differential in the total complex, carry information about all elementary physical processes in this system. 

2. If we cut the tree along red edges, we get a decomposition of our tree T into several trees. According to the theorem of O. Lyashko (\cite{lya}), such decomposition, if T is a Dynkin diagram of a simple singularity, corresponds to the decay of the initial singularity into several ones, whose Dynkin diagrams correspond to the parts of the tree T.

Therefore the differential $D$ may be interpreted as follows. Let $W$ be a space with the basis $\{v_{+},v_{-}\}$
and fixed skew form $\eta=v_{+}\wedge v_{-}$. Now to each decomposition of tree T we assign a vector space: it is $W$
in the tensor power equal to the number of free vertices. If we join two parts, we delete one red edge, and hence create two new vertices. Now the operator $D$ corresponds to the map $m$ from the tensor product of two vector spaces on parts to the vector space on their union.  We define
$$m(a,b)=a\otimes \eta\otimes b.$$
This gives a natural possibility of generalization of the construction of Seifert cohomology by changing $W$ by some other space.

3. An interesting question is to define the cohomology theory that categorifies the Alexander polynomial for general graphs, not only for trees. To be more precise, it would be interesting to have a theory independent of the choice of a basis:  it is well known (e. g. \cite{book}) that different distinguished basises in the vanishing homology of a singularity are related by the braid group action. One can compute the effect of braid group action on the Dynkin diagram, and one can ask for a cohomology theory that gives the same answers for different Dynkin diagrams of the same singularity in different basises. Even for simple singularities one can change the basis such that the Dynkin diagram will not be a tree.

Since the Heegard-Floer homology are independent of such choices, a natural question for such a conjectural theory will be about its relation to the Heegard-Floer theory.

4. In \cite{tl} V. Toledano-Laredo constructed a homology theory, also starting  from the Dynkin diagram. We do not know any relation between his homology and ours despite that the size of his complex looks similar to the size of ours one.

Moscow State University, Department of Mathematics and Mechanics.

E. mail: gorsky@mccme.ru

\end{document}